\documentclass[11pt]{article}

\usepackage[margin=1in]{geometry}
\usepackage{color}
\usepackage{amssymb}
\usepackage{amsmath}
\usepackage{amsthm}
\usepackage{url}
\usepackage{framed}
\usepackage{enumitem}
\usepackage{graphicx}
\usepackage{multicol}
\usepackage{longtable}

\usepackage[T1]{fontenc}
\usepackage[utf8]{inputenc}
\usepackage{authblk}


\newtheorem{theorem}{Theorem}

\newtheorem{conjecture}[theorem]{Conjecture}

\begin{document}

\title{Tables of the existence of equiangular tight frames}

\author{Matthew Fickus \quad and \quad Dustin G.\ Mixon\thanks{Email: dustin.mixon@gmail.com}}
\affil{Department of Mathematics and Statistics, Air Force Institute of Technology, Wright-Patterson AFB, OH 45433}

\maketitle

\begin{abstract}
A Grassmannian frame is a collection of unit vectors which are optimally incoherent.
To date, the vast majority of explicit Grassmannian frames are equiangular tight frames (ETFs).
This paper surveys every known construction of ETFs and tabulates existence for sufficiently small dimensions.
\end{abstract}

\section{Introduction}

Let $\Phi=\{\varphi_n\}_{n=1}^N$ denote a finite sequence of points from the unit sphere in $\mathbb{R}^M$ or $\mathbb{C}^M$.
We want $\Phi$ to minimize the \textbf{worst-case coherence}, defined as
\[
\mu(\Phi)
:=\max_{\substack{i,j\in\{1,\ldots,N\}\\i\neq j}}|\langle \varphi_i,\varphi_j\rangle|.
\]
Minimizers of worst-case coherence are called \textbf{Grassmannian frames}~\cite{StrohmerH:03}.
A compactness argument establishes that a Grassmannian frame exists for every pair of parameters $(M,N)$.
But how might one construct a Grassmannian frame?

The \textbf{Welch bound}~\cite{Welch:74,StrohmerH:03} is a lower bound on the worst-case coherence:
\[
\mu(\Phi)
\geq\sqrt{\frac{N-M}{M(N-1)}}.
\]
The proof of the Welch bound is instructive:
Identifying $\Phi$ with the $M\times N$ matrix $[\varphi_1\cdots\varphi_N]$ gives
\[
N+N(N-1)\mu(\Phi)^2
\geq\sum_{i=1}^N\sum_{j=1}^N|\langle\varphi_i,\varphi_j\rangle|^2
=\|\Phi^*\Phi\|_\mathrm{HS}^2
\geq\frac{N^2}{M},
\]
where the last inequality follows from rearranging $0\leq\|\Phi\Phi^*-\frac{N}{M}I\|_\mathrm{HS}^2=\|\Phi^*\Phi\|_\mathrm{HS}^2-\frac{N^2}{M}$.
As such, equality in the Welch bound occurs precisely when
\begin{itemize}
\item[(i)] the cosines $|\langle\varphi_i,\varphi_j\rangle|$ are the same for every $i$ and $j\neq i$, and
\item[(ii)] the frame $\Phi$ is $\textbf{tight}$, that is $\Phi\Phi^*=\frac{N}{M}I$.
\end{itemize}
Such ensembles are called \textbf{equiangular tight frames}, and since they achieve equality in a lower bound of $\mu(\Phi)$, we conclude that they are necessarily Grassmannian.
This result is significant because Grassmannian frames are very difficult to identify in general (see~\cite{BenedettoK:06}, for example\footnote{Also, see \cite{BodmannH:15,CasazzaH:16} for recent developments in non-ETF Grassmannian frames.}), and the additional structural information afforded by ETFs make them more accessible.
Indeed, as we will see, there are currently several infinite families of known ETFs.

Examples of ETFs include the cube roots of unity (viewed as vectors in $\mathbb{R}^2$) and the vertices of the origin-centered tetrahedron (viewed as vectors in $\mathbb{R}^3$).
The apparent beauty of ETFs coupled with their importance as Grassmannian frames has made them the subject of active research recently.
To date, several tables of the existence of ETFs can be found in the literature, e.g., \cite{SustikTDH:07,Waldron:09,Redmond:09}.
Instead of publishing a fixed survey and table of known ETFs, the present document provides a ``living'' alternative so as to account for recent and future developments in the construction or impossibility of ETFs in various dimensions.
The intent is to update this document periodically as developments arrive.

The following section summarizes the known infinite families of ETFs.
Sections~3 and~4 then provide more detail in the cases of real and complex ETFs, respectively, and Section~5 focuses on ETFs with redundancy close to 2.
We discuss our methodology for constructing existence tables in Section~6, and the tables can be found at the end of the paper.

\section{Infinite families of equiangular tight frames}

What follows is a list of trivial ETFs (we will ignore these for the rest of the paper):
\begin{itemize}
\item
\textbf{Orthonormal bases.}
This case takes $N=M$, and it is easy to verify the ETF conditions.
\item
\textbf{Regular simplices.}
In this case, $N=M+1$.
For a simple construction of this example, take $N-1$ rows from an $N\times N$ discrete Fourier transform matrix.
Then the resulting columns, after being scaled to have unit norm, form an ETF.
\item
\textbf{Frames in one dimension.}
When $M=1$, any unit norm frame amounts to a list of scalars of unit modulus, and such frames are necessarily ETFs.
\end{itemize}

To date, there are only a few known infinite families of nontrivial ETFs.
Interestingly, despite ETFs being characterized in terms of functional analysis (namely, equality in the Welch bound), each of the known infinite families is based on some sort of combinatorial design:
\begin{itemize}
\item
\textbf{Strongly regular graphs.}
A $(v,k,\lambda,\mu$)-strongly regular graph is a $k$-regular graph with $v$ vertices such that every pair of adjacent vertices has $\lambda$ common neighbors, whereas every pair of non-adjacent vertices has $\mu$ common neighbors.
One may manipulate the adjacency matrix of a strongly regular graph with appropriate parameters to find an embedding of some $M$-dimensional real ETF of $N=v+1$ vectors in $\mathbb{R}^N$ (here, $M$ is a complicated function of the graph parameters; see~\cite{Waldron:09} for details).
In fact, real ETFs are in one-to-one correspondence with a subclass of strongly regular graphs in this way.
As an example, the case where $N=2M$ corresponds to the so-called conference graphs.
This graph-based construction has since been generalized in multiple ways to produce complex ETFs (namely, using antisymmetric conference matrices~\cite{Renes:07,Strohmer:08}, complex Hadamard matrices~\cite{Szollosi:13}, or distance-regular antipodal covers of the complete graph~\cite{CountinhoGSZ:16,FickusJMPW:16}).
\item
\textbf{Difference sets.}
Let $G$ be a finite abelian group.
Then $D\subseteq G$ is said to be a $(G,k,\lambda)$-difference set if $|D|=k$ and for every nonzero $g\in G$, there are exactly $\lambda$ different pairs $(d_1,d_2)\in D\times D$ such that $g=d_1-d_2$.
One may use any difference set to construct an ETF with $M=|D|$ and $N=|G|$ (see~\cite{StrohmerH:03,XiaZG:05,DingF:07}).
In particular, each vector in the ETF is obtained by taking a character of $G$ and restricting its domain to $D$ (before scaling to have unit norm).
\item
\textbf{Steiner systems.}
A $(2,k,v)$-Steiner system is a $v$-element set $S$ of points together with a collection $\mathcal{B}$ of $k$-element subsets of $S$ called blocks with the property that each $2$-element subset of $S$ is contained in exactly one block.
It is not difficult to show that each point is necessarily contained in exactly $r=(v-1)/(k-1)$ blocks.
One may use any $(2,k,v)$-Steiner system to construct an ETF in $\mathbb{C}^\mathcal{B}$ (see~\cite{FickusMT:12}).
Specifically, for each point $p\in S$, embed an $r$-dimensional regular simplex into $\mathbb{C}^\mathcal{B}$ so as to be supported on the blocks that contain $p$.
The union of these embedded simplices then forms an ETF of $v(r+1)$ vectors in $\mathbb{C}^\mathcal{B}$.
Every such construction necessarily has $N>2M$.
This construction has been modified in two different ways:
\cite{FickusJMP:16} uses Steiner ETFs from triples systems to generalize Example~7.10 in Tremain's notes~\cite{Tremain:09}, and \cite{FickusMJ:16} provides other Steiner-like constructions involving projective planes that contain a hyperoval.
A generalization of Steiner systems, namely, quasi-symmetric designs, have also been used to construct ETFs~\cite{FickusJMP:15}.
\end{itemize}

\section{Real equiangular tight frames}

In this section, we describe what is known about real equiangular tight frames.
Throughout, we use $\exists\operatorname{RETF}(M,N)$ to denote the statement ``there exists a real equiangular tight frame with parameters $(M,N)$.''
We start with some basic properties:

\begin{theorem}[see~\cite{SustikTDH:07}]
\label{thm.basic bounds}
$\exists\operatorname{RETF}(M,N)$ implies each of the following:
\begin{itemize}
\item[(a)]
$N\leq M(M+1)/2$.
\item[(b)]
$\exists\operatorname{RETF}(N-M,N)$.
\end{itemize}
\end{theorem}

Part~(a) above can be seen by observing that the rank-$1$ matrices $\{\varphi_n\varphi_n^*\}_{n=1}^N$ are necessarily linearly independent in the $M(M+1)/2$-dimensional space of $M\times M$ symmetric matrices (this follows from computing the spectrum of their Gram matrix).
Part~(b) uses a concept known as the \textbf{Naimark complement}.
In particular, since any ETF $\Phi$ satisfies $\Phi\Phi^*=\frac{N}{M}I$, the rows of $\Phi$ can be viewed as orthonormal vectors (suitably scaled).
As such, one may complete the orthonormal basis with $N-M$ other row vectors.
Collecting these rows into an $(N-M)\times N$ matrix and normalizing the columns results in another ETF (called the Naimark complement of the original).

Real ETFs are intimately related to graphs.
Given any real ETF, negate some of the vectors so that each one has positive inner product with the last vector (this process produces another ETF).
Next, remove the last vector to get a subcollection of vectors $\Psi$ (this is no longer an ETF).
Use $\Psi=\{\psi_n\}_{n=1}^{N-1}$ to build a graph in the following way:
Take $v=N-1$ vertices and say vertex $i$ is adjacent to vertex $j$ if $\langle \psi_i,\psi_j\rangle<0$.
It turns out that this graph is necessarily strongly regular with parameters determined by $M$ and $N$:

\begin{theorem}[Corollary~5.6 in~\cite{Waldron:09}]
\label{thm.real etfs are srgs}
$\exists\operatorname{RETF}(M,N)$ if and only if $\exists\operatorname{SRG}(N-1,k,\frac{3k-N}{2},\frac{k}{2})$ with
\[
k=\frac{N}{2}-1+\bigg(1-\frac{N}{2M}\bigg)\sqrt{\frac{M(N-1)}{N-M}}.
\]
\end{theorem}

The spectrum of a strongly regular graph can be expressed in terms of its graph parameters.
In fact, it turns out that the eigenvalues must be integer, which in turn implies the following integrality conditions:

\begin{theorem}[Theorem~A in~\cite{SustikTDH:07}]
\label{thm.integrality for redundancy not 2}
Suppose $N\neq2M$.
Then $\exists\operatorname{RETF}(M,N)$ implies that
\[
\sqrt{\frac{M(N-1)}{N-M}},
\qquad
\sqrt{\frac{(N-M)(N-1)}{M}}
\]
are both odd integers.
\end{theorem}

Since we identify real ETFs with certain strongly regular graphs, we can leverage necessary conditions for existence of the latter to inform existence of the former:

\begin{theorem}[see~\cite{Brouwer:07,BrouwerH:12}]
\label{thm.graph conditions}
Given $v$, $k$, $\lambda$ and $\mu$, let $r\geq0$ and $s\leq-1$ denote the solutions to
\[
x^2+(\mu-\lambda)x+(\mu-k)=0,
\]
and take
\[
f:=\frac{s(v-1)+k}{s-r},
\qquad
q_{11}^1:=\frac{f^2}{v}\bigg(1+\frac{r^3}{k^2}-\frac{(r+1)^3}{(v-k-1)^2}\bigg).
\]
Then $\exists\operatorname{SRG}(v,k,\lambda,\mu)$ implies each of the following:
\begin{itemize}
\item[(a)]
The \textbf{Krein conditions} are satisfied:
\[
(r+1)(k+r+2rs)\leq(k+r)(s+1)^2,
\qquad
(s+1)(k+s+2rs)\leq(k+s)(r+1)^2.
\]
\item[(b)]
$\displaystyle{v\leq\left\{\begin{array}{ll}\frac{1}{2}f(f+3)&\mbox{if }q_{11}^1=0\\\frac{1}{2}f(f+1)&\mbox{if }q_{11}^1\neq0\end{array}\right.}$.
\item[(c)]
If $\mu=1$, then $\frac{vk}{(\lambda+1)(\lambda+2)}$ is integer.
\end{itemize}
\end{theorem}

In the case of real ETFs, $\mu=k/2$, and so $\mu=1$ implies $k=2$, thereby implying $3-N/2=\lambda\geq0$, i.e., $N\leq 6$.
However, an exhaustive search shows that $\mu>1$ for every $(M,N)$ satisfying the conditions in Theorems~\ref{thm.basic bounds} and~\ref{thm.integrality for redundancy not 2} with $N\leq 6$, and so part (c) above is moot.
It is unclear whether part (a) or part (b) is also covered by the previous conditions; indeed, for the entire range of $(M,N)$ pairs we tested for our tables, these necessary conditions failed to rule out any pairs $(M,N)$ made plausible by the previous conditions.
Interestingly, for part (b), the case $q_{11}^1=0$ is important when discerning the existence of real ETFs; for example, when $(M,N)=(21,28)$ and $(253,276)$, $v$ lies between $\frac{1}{2}f(f+1)$ and $\frac{1}{2}f(f+3)$, but $q_{11}^1=0$ in both instances (also, real ETFs are known to exist in both instances).

\subsection{Maximal real ETFs}

Let $S^{M-1}$ denote the unit sphere in $\mathbb{R}^M$.
A \textbf{spherical $t$-design} in $S^{M-1}$ is a finite collection $X\subseteq S^{M-1}$ satisfying
\[
\frac{1}{|X|}\sum_{x\in X}f(x)
=\frac{1}{\omega_{M-1}}\int_{S^{M-1}}f(x)dS
\]
for every polynomial $f(x)$ of degree at most $t$; on the right-hand side, the integral is taken with respect to the Haar measure of $S^{M-1}$, and $\omega_{M-1}$ denotes the measure of $S^{M-1}$.
For each $t$, there is a Fisher-type inequality that provides a lower bound on $|X|$.
For example, a spherical $5$-design $X$ necessarily satisfies $|X|\geq M(M+1)$~\cite{GoethalsS:79}.
We say a spherical $t$-design is \textbf{tight} if it satisfies equality in the Fisher-type inequality.

\begin{theorem}[see \cite{Redmond:09,Mixon:16}]
\label{thm.real maximal etfs are spherical 5-designs}
Every tight spherical $5$-design $X$ is of the form $X=\Phi\cup(-\Phi)$ for some ETF $\Phi$ of $M(M-1)/2$ elements in $\mathbb{R}^M$.
Conversely, for every such ETF $\Phi$, the collection $\Phi\cup(-\Phi)$ forms a tight spherical $5$-design.
\end{theorem}

In this special case where $N=M(M+1)/2$, it is straightforward to verify that Theorem~\ref{thm.integrality for redundancy not 2} implies something special about the form of $M$.
In particular, provided $N\neq2M$ (i.e., $M\neq3$), $\exists\operatorname{RETF}(M,M(M+1)/2)$ requires an integer $m\geq1$ such that $M=(2m+1)^2-2$.
Overall, $\exists\operatorname{RETF}(M,M(M+1)/2)$ implies $M\in\{3,7,23,47,79,119,167,\ldots\}$.
The following theorem summarizes what is known about the existence of these ETFs:

\begin{theorem}[see~\cite{GoethalsS:75,Makhnev:02,BannaiMV:05}]
\label{thm.maximal real etfs}
\
\begin{itemize}
\item[(a)]
$M\in\{3,7,23\}$ implies $\exists\operatorname{RETF}(M,M(M+1)/2)$.
\item[(b)]
$M=47$ implies $\nexists\operatorname{RETF}(M,M(M+1)/2)$.
\item[(c)]
Suppose $k\equiv2\bmod3$, $k$ and $2k+1$ are both square-free, and take $m=2k$. 
Then $M=(2m+1)^2-2$ implies $\nexists\operatorname{RETF}(M,M(M+1)/2)$.
\end{itemize}
\end{theorem}

Part~(b) above was originally proved by Makhnev~\cite{Makhnev:02} in terms of strongly regular graphs, and soon thereafter, Bannai, Munemasa and Venkov~\cite{BannaiMV:05} found an alternative proof in terms of spherical $5$-designs (along with a proof of part~(c) above).
Other than the dimension bound (Theorem~\ref{thm.basic bounds}) and the integrality conditions (Theorem~\ref{thm.integrality for redundancy not 2}), this was the first nonexistence result\footnote{In the time since, it has been proven that $\nexists\operatorname{RETF}(19,76)$ and $\nexists\operatorname{RETF}(20,96)$; see~\cite{AzarijaM:15,AzarijaM:16}.} for real ETFs with $N>2M$.
In fact, this disproved a conjecture that was posed in~\cite{SustikTDH:07} and reiterated in~\cite{Waldron:09}.

\section{Complex equiangular tight frames}

Far less is known about complex equiangular tight frames.
In this section, we use $\exists\operatorname{ETF}(M,N)$ to denote the statement ``there exists an equiangular tight frame with parameters $(M,N)$.''
The following result is analogous to Theorem~\ref{thm.basic bounds} (as is the proof):

\begin{theorem}[see~\cite{SustikTDH:07}]
\label{thm.basic bounds for complex}
$\exists\operatorname{ETF}(M,N)$ implies each of the following:
\begin{itemize}
\item[(a)]
$N\leq M^2$.
\item[(b)]
$\exists\operatorname{ETF}(N-M,N)$.
\end{itemize}
\end{theorem}

Theorem~\ref{thm.basic bounds for complex} is the only known general necessary condition on the existence of ETFs (unlike the real case, which enjoys conditions like Theorem~\ref{thm.integrality for redundancy not 2}).
In fact, there is only one case in which nonexistence has been proved beyond Theorem~\ref{thm.basic bounds for complex}:
There does not exist an ETF with $(M,N)=(3,8)$ (nor with $(5,8)$ by the Naimark complement)~\cite{Szollosi:14}.
Unfortunately, Sz\"{o}ll\H{o}si's proof of this fact involves the computation of a Gr\"{o}bner basis, which doesn't seem to provide a more general nonexistence result.

In the absence of techniques for nonexistence results, the dominant mode of progress in the complex case has been existence results by various constructions.
The remainder of this section details some of these constructions.

\subsection{Difference sets}
\label{subsection.difference sets}

As discussed in Section~2, a difference set in an abelian group immediately yields an ETF.
Thankfully, there is an entire community of researchers who have developed a rich theory of difference sets (see~\cite{JungnickelPS:07} for an overview).
The following summarizes most of the existence results:

\begin{theorem}[see~\cite{JungnickelPS:07}]
\label{thm.difference sets make etfs}
Each of the following implies $\exists\operatorname{ETF}(M,N)$:
\begin{itemize}
\item[(a)]
$\operatorname{PG}(m-1,q)$: $M=\frac{q^{m-1}-1}{q-1}$, $N=\frac{q^m-1}{q-1}$ for some prime power $q$ and integer $m\geq3$.
\item[(b)]
$\operatorname{MF}(q,d)$: $M=q^d\cdot\frac{q^{d+1}-1}{q-1}$, $N=q^{d+1}(1+\frac{q^{d+1}-1}{q-1})$ for some prime power $q$ and integer $d\geq1$.
\item[(c)]
$\operatorname{Sp}(d)$: $M=3^d\cdot\frac{3^{d+1}+1}{2}$, $N=3^{d+1}\cdot\frac{3^{d+1}-1}{2}$ for some integer $d\geq1$.
\item[(d)]
$\operatorname{Paley}(t)$: $M=2t-1$, $N=4t-1$, provided $N$ is a prime power.
\item[(e)]
$\operatorname{Cyclo}_{4}(t)$: $M=t^2$, $N=4t^2+1$ for some odd integer $t$, provided $N$ is a prime power.
\item[(f)]
$\operatorname{Cyclo}_{4'}(t)$: $M=t^2+3$, $N=4t^2+9$ for some odd integer $t$, provided $N$ is a prime power.
\item[(g)]
$\operatorname{Cyclo}_{8}(t,u)$: $M=t^2=8u^2+1$, $N=8t^2+1=64u^2+9$ for some odd integers $t$ and $u$, provided $N$ is a prime power.
\item[(h)]
$\operatorname{Cyclo}_{8'}(t,u)$: $M=t^2+7=8u^2+56$, $N=8t^2+49=64u^2+441$ for some odd integer $t$ and even integer $u$, provided $N$ is a prime power.
\item[(i)]
$\operatorname{H}(t)$: $M=2t^2+13$, $N=4t^2+27$, provided $N\equiv1\bmod6$ is a prime power.
\item[(j)]
$\operatorname{TPP}(q)$: $M=\frac{q^2+2q-1}{2}$, $N=q^2+2q$, provided $q$ and $q+2$ are both odd prime powers.
\end{itemize}
\end{theorem}

We note that other infinite families of difference sets include the Davis--Jedwab--Chen and Hadamard difference sets; the conditions for the existence of these difference sets are much more complicated, and so we omit them.
We also note an apparent typographical error in~\cite{JungnickelPS:07}:
Hall difference sets (which form the basis of part~(i) above) should be defined in terms of powers of a primitive element of a finite field instead of powers of arbitrary units.
Additionally, the twin prime powers $q$, $q+2$ used in part~(j) above need to be odd (as stated above)~\cite{BethJL:99}, though this condition is missing in~\cite{JungnickelPS:07}.

\subsection{Maximal ETFs}

Recall from the real case that ETFs with the maximum possible number of vectors found applications in cubature (Theorem~\ref{thm.real maximal etfs are spherical 5-designs}), though they don't always exist (Theorem~\ref{thm.maximal real etfs}).
Considering Theorem~\ref{thm.basic bounds for complex}, the analogous ETFs in the complex case would have $N=M^2$, and in fact, these particular ETFs also enjoy special applications.
In particular, in quantum mechanics, these are called \textbf{symmetric, informationally complete, positive operator--valued measures (SIC-POVMs)}, and in this context, they are foundational to the theory of quantum Bayesianism~\cite{FuchsS:11} and find further applications in both quantum state tomography~\cite{CavesFS:02} and quantum cryptography~\cite{FuchsS:03}.
Needless to say, there is considerable interest in the existence of maximal ETFs, and (unlike the real case) they are conjectured to exist in every dimension:

\begin{conjecture}[Zauner's conjecture~\cite{Zauner:99}]
$\exists\operatorname{ETF}(M,M^2)$ for every $M\geq2$.
\end{conjecture}

Unfortunately, there is currently no infinite family of known maximal ETFs.
The following summarizes what is known:

\begin{theorem}[see~\cite{Zauner:online,Chien:14}]
\label{thm.maximal etfs}
If $M\leq 17$ or $M\in\{19,24,28,35,48\}$, then $\exists\operatorname{ETF}(M,M^2)$.
\end{theorem}

For the record, Scott and Grassl~\cite{ScottG:10} ran numerical tests to find ensembles of $M^2$ vectors in $\mathbb{C}^M$ that are within machine precision of satisfying the ETF conditions, thereby suggesting that Zauner's conjecture is likely true up to dimension 67 (at least).

\section{Redundancy 2, more or less}

In this section, we focus on the special case of real and complex ETFs whose redundancy $N/M$ is close to 2.
First, note that Theorem~\ref{thm.integrality for redundancy not 2} does not apply when $N=2M$.
In this case, Theorem~\ref{thm.real etfs are srgs} gives that a real ETF corresponds to a strongly regular graph with $k=\frac{N}{2}-1$, in which case the graph is called a \textbf{conference graph}.
A bit more is known about this special type of strongly regular graph:

\begin{theorem}[see~\cite{Belevitch:50}]
\label{thm.red 2 integrality}
$\exists\operatorname{RETF}(M,2M)$ implies that $M$ is odd and $2M-1$ is a sum of two squares.
\end{theorem}

A \textbf{conference matrix} is an $N\times N$ matrix $C$ with zero diagonal and off-diagonal entries $\pm1$ such that $CC^\top=(N-1)I$.
Given a conference graph with vertices $\{1,\ldots,N-1\}$, one may construct a conference matrix $C$ with row and column indices in $\{0,1,\ldots,N-1\}$ such that $C_{ii}=0$ for every $i$, $C_{ij}=-1$ whenever $i$ and $j$ are nonzero and adjacent in the graph, and $C_{ij}=1$ otherwise.
Then $C$ gives the sign pattern of the Gram matrix of the ETF which corresponds to the conference graph.
In fact, since the rows and columns of any symmetric conference matrix can be signed to have this form, then by Theorem~\ref{thm.real etfs are srgs}, we identify all symmetric conference matrices with the Gram matrices of real ETFs with redundancy 2.

In the case of complex ETFs, we want the sign pattern of the Gram matrix to be self-adjoint instead of symmetric.
To accomplish this, one might multiply an antisymmetric conference matrix by $\mathrm{i}=\sqrt{-1}$ (indeed, this leads to a complex ETF, as one might expect).
The following theorem summarizes the parameters of known symmetric and antisymmetric conference matrices, which in turn lead to ETFs of redundancy 2:

\begin{theorem}[see~\cite{IoninK:07,BaloninS:14}]
\label{thm.conference constructions}
Each of the following implies $\exists\operatorname{RETF}(M,N)$ with $N=2M$:
\begin{itemize}
\item[(a)]
$N=q+1$ for some prime power $q\equiv1\bmod4$.
\item[(b)]
$N=q^2(q+2)+1$ for some prime power $q\equiv3\bmod4$ and prime power $q+2$.
\item[(c)]
$N=5\cdot9^{2t+1}+1$ for some integer $t\geq0$.
\item[(d)]
$N=(h-1)^{2s}+1$ for some integer $s\geq1$, where $h$ is the order of a skew-Hadamard matrix.
\item[(e)]
$N=(n-1)^s+1$ for some integer $s\geq2$, where $n$ is the order of a conference matrix.
\end{itemize}
Each of the following implies $\exists\operatorname{ETF}(M,N)$ with $N=2M$:
\begin{itemize}
\item[(f)]
$N=q+1$ for some prime power $q\equiv3\bmod4$.
\item[(g)]
$N=h$, where $h\geq4$ is the order of a skew-Hadamard matrix.
\end{itemize}
\end{theorem}

Interestingly, in the special case where the ETF corresponds to an antisymmetric conference matrix (meaning all the inner products are $\pm\mathrm{i}$ times the Welch bound), one can perform the following operation to produce another ETF~\cite{Renes:07,Strohmer:08}:
Remove any one of the vectors to get a $M\times(2M-1)$ matrix $\Psi$, and then compute $\Phi=(\alpha\Psi\Psi^*)^{-1/2}\Psi$, where $\alpha=M/(2M-1)$; then the columns of $\Phi$ form an ETF.
As such, statements (f) and (g) in Theorem~\ref{thm.conference constructions} also imply $\exists\operatorname{ETF}(N/2,N-1)$.

For another variation on ETFs of redundancy 2, recall that a \textbf{complex Hadamard matrix} is an $N\times N$ matrix $H$ with entries of unit modulus satisfying $HH^*=NI$.
Complex Hadamard matrices are related to conference matrices:
If $C$ is symmetric, then $\mathrm{i}I+C$ is complex Hadamard, and if $C$ is antisymmetric, then $I+C$ (and therefore $\mathrm{i}I+\mathrm{i}C$) is complex Hadamard.
In general, when $N/M$ is sufficiently close to 2 (explicitly, when $M$ is between $(N-\sqrt{N})/2$ and $(N+\sqrt{N})/2$), then $I+\mu Q$ is the Gram matrix of an $M\times N$ ETF (with $\mu$ denoting the Welch bound) if and only if $\lambda I+Q$ is a complex Hadamard matrix for some $\lambda$~\cite{Szollosi:13}.
As such, we identify ETFs in this range with complex Hadamard matrices of constant diagonal and self-adjoint off-diagonal. 

Beyond Theorem~\ref{thm.conference constructions}, there are two known families of complex Hadamard matrices with these specifications.
First, in the case where the Hadamard matrix is real, we want a symmetric Hadamard matrix with constant diagonal (these are known as \textbf{graphical Hadamard matrices}).
Such a matrix $H$ has positive diagonal (without loss of generality), and so $\operatorname{Tr}(H)=N$ and $H^2=HH^\top=NI$, which together imply that the eigenvalues are $\pm\sqrt{N}$ with multiplicities $(N\pm\sqrt{N})/2$, respectively.
As such, putting $Q=H-I$, then the Gram matrix $I+\mu Q$ will have rank $M=(N+\sqrt{N})/2$.
A lot of work in the Hadamard literature has focused on a special class of these matrices, namely, \textbf{regular symmetric Hadamard matrices of constant diagonal (RSHCD)}, where ``regular'' indicates that the all-ones vector is an eigenvector.
The following result summarizes the status of this literature:

\begin{theorem}[see~\cite{BrouwerH:12}]
\label{thm.RSHCD}
Each of the following implies $\exists\operatorname{RETF}(M,N)$ with $M=(N+\sqrt{N})/2$:
\begin{itemize}
\item[(a)]
$N\in\{4,36,100,196\}$.
\item[(b)]
$N=h^2$, where $h$ is the order of a Hadamard matrix.
\item[(c)]
$N=a^2$, where $a-1$ and $a+1$ are both odd prime powers.
\item[(d)]
$N=a^2$, where $a+1$ is a prime power and $a$ is the order of a symmetric conference matrix.
\item[(e)]
$N=4t^2$, provided there exists a set of $t-2$ mutually orthogonal Latin squares of order $2t$.
\item[(f)]
$N=4t^4$, for some $t\geq1$.
\end{itemize}
\end{theorem}

As indicated above, there is a second infinite family of complex Hadamard matrices with constant diagonal and self-adjoint off-diagonal, namely \textbf{self-adjoint complex Hadamard matrices with constant diagonal (SCHCD)}.
(Implicitly, the diagonal of such matrices is real.)
In particular, \cite{Szollosi:13} constructs SCHCDs of order $n^2$ for every $n\geq2$.

\section{Table methodology}

The tables that annotate the existence of ETFs can be found at the end of this paper.
In this section, we describe how these tables were generated.

\subsection{Table~\ref{table.real etfs red>2}}
\label{subsection.how to real}

Here, $\alpha$ denotes the inverse coherence, namely $\sqrt{M(N-1)/(N-M)}$, and $k$ denotes the strongly regular graph parameter from Theorem~\ref{thm.real etfs are srgs}.
To build this table, we used the following procedure:
\begin{enumerate}
\item
\textbf{Find parameters that satisfy integrality.}
We first found all pairs $(M,N)$ with $M/2<N\leq M(M+1)/2$ and $N\leq1300$ that satisfy the integrality conditions (Theorem~\ref{thm.integrality for redundancy not 2}).
Indeed, we chose to treat the case where $N=2M$ in a different table since it uses fundamentally different necessary conditions.
In light of the Naimark complement, we also focused on $N>2M$ (see Theorem~\ref{thm.basic bounds}).
The upper bound $N\leq M(M+1)/2$ also follows from Theorem~\ref{thm.basic bounds}, and the bound $N\leq1300$ was a somewhat arbitrary choice, but it coincides with the table of strongly regular graphs found in~\cite{Brouwer:15}.
\item
\textbf{Test necessary conditions for strongly regular graphs.}
Next, we ensure that the strongly regular graph parameters corresponding to both $(M,N)$ and $(N-M,N)$ satisfy the conditions in Theorem~\ref{thm.graph conditions}.
Any violation would be reflected in the Notes column of the table as ``failed graph test,'' but no such violation occurred.
\item
\textbf{Find parameters of real ETFs arising from difference sets.}
Considering Corollary~2 in~\cite{JasperMF:14}, a difference set can be used to produce a real ETF only if $M=2^j(2^{j+1}\pm1)$ and $N=2^{2j+2}$.
Conversely, as explained in~\cite{DingF:07}, McFarland difference sets can be used to construct a real ETF for every such $(M,N)$.
\item
\textbf{Find parameters of real ETFs arising from Steiner systems.}
A few infinite families of Steiner systems are described in~\cite{FickusMT:12} to construct ETFs.
Of these, real ETFs can be constructed with Steiner systems of 2-blocks, 3-blocks, 5-blocks or 6-blocks\footnote{In order to consider 6-blocks, we needed to hard-code certain exceptions and unknown cases detailed in~\cite{AbelG:07}.}, as well as those constructed from affine and projective geometries.
We did not consider any other Steiner systems.
\item
\textbf{Find parameters of real ETFs arising from other constructions.}
Here, we account for constructions arising from generalized quadrangles~\cite{FickusJMPW:16}, RSHCDs (see Theorem~\ref{thm.RSHCD}), Tremain ETFs~\cite{FickusJMP:16}, and quasi-symmetric designs~\cite{FickusJMP:15}. When implementing Theorem~\ref{thm.RSHCD}, we did not consider (e), as this is difficult to test.
\item
\textbf{Find parameters for which real ETFs do not exist.}
Next, we applied Theorem~\ref{thm.maximal real etfs} along with~\cite{AzarijaM:15,AzarijaM:16} to tabulate the nonexistence of certain ETFs, denoted by ``DNE.''
\item
\textbf{Find parameters of real ETFs arising from strongly regular graphs.}
At this point, the entire ETF literature has been exhausted, and so we appeal to Theorem~\ref{thm.real etfs are srgs} and the table of strongly regular graphs found in~\cite{Brouwer:15} to fill in any remaining known constructions.
This revealed four remaining constructions, which we denote in the table by ``SRG.''
\end{enumerate}

\subsection{Table~\ref{table.etfs red>2}}
\label{subsection.how to complex}

In the complex case, there are no known integrality conditions that one can use to narrow down the list of possible dimensions.
In fact, beyond the standard dimension bounds, there is only one result on the nonexistence of complex ETFs:
There does not exist an ETF with $(M,N)=(3,8)$ or $(5,8)$~\cite{Szollosi:14}.
This feature of the complex case informs how one ought to build a table.
Indeed, Table~\ref{table.etfs red>2} only provides parameters for which ETFs are known to exist.

As in the real case, this table only gives the known ETFs for which $N>2M$ (with the understanding that ETFs with parameters $(N-M,N)$ also necessarily exist by the Naimark complement).
In addition, we only considered $M\leq300$ and $N\leq1300$.
These restrictions actually come from existing tables of difference sets~\cite{Gordon:online} and strongly regular graphs~\cite{Brouwer:15}, respectively, meaning Table~\ref{table.etfs red>2} accounts for every known complex ETF in this range.

The following describes the methodology we used to construct Table~\ref{table.etfs red>2}:
\begin{enumerate}
\item
\textbf{Find parameters of ETFs arising from difference sets.}
We first collected ETF parameters from each part of Theorem~\ref{thm.difference sets make etfs}.
By doing this, we specifically neglected both Davis--Jedwab--Chen and Hadamard difference sets, as these are particularly difficult to parameterize.
However, since our table only takes $M\leq300$, we were able to consult the La Jolla Difference Set Repository~\cite{Gordon:online} and identify any difference sets in this range that we missed.
This process revealed 6 additional $(M,N)$ pairs that we then hard-coded into the table.
\item
\textbf{Find parameters of ETFs arising from Steiner systems.}
Similar to our process for the real ETF tables, we appealed to~\cite{FickusMT:12} to identify infinite families of ETFs which arise from Steiner systems (as the table illustrates, many more Steiner systems yield complex ETFs).
\item
\textbf{Find parameters of ETFs arising from other constructions.}
Here, we account for constructions arising from modifying skew-symmetric conference matrices~\cite{Renes:07,Strohmer:08}, generalized quadrangles~\cite{FickusJMPW:16}, RSHCDs (see Theorem~\ref{thm.RSHCD}), SCHCDs~\cite{Szollosi:13}, Tremain ETFs~\cite{FickusJMP:16}, quasi-symmetric designs~\cite{FickusJMP:15}, and hyperovals~\cite{FickusMJ:16}.
We only considered skew-symmetric conference matrices indicated by (f) and (g) in Theorem~\ref{thm.conference constructions}, and for (g), we only considered skew-Hadamard matrices of order $h=2^t$ for some $t\geq2$.
\item
\textbf{Find parameters of maximal ETFs.}
As indicated in Theorem~\ref{thm.maximal etfs}, there are a few cases in which maximal ETFs are known to exist, and we hard-coded these into the tables.
\item
\textbf{Find parameters of real ETFs arising from strongly regular graphs.}
Finally, like Table~\ref{table.real etfs red>2}, we considered the remaining ETFs that come from the table of strongly regular graphs~\cite{Brouwer:15}.
\end{enumerate}

\subsection{Table~\ref{table.etfs red=2}}

We followed a similar process to build a table for the special case where $N=2M$, namely Table~\ref{table.etfs red=2}.
In this case, we first determined which $M\leq 150$ satisfy Theorem~\ref{thm.red 2 integrality}.
Next, we found the $N$ which satisfy each part of Theorem~\ref{thm.conference constructions}.
The column labeled ``$\mathbb{R}$?'' gives a ``+'' if there exists a known real construction, ``-'' if the parameters fail to satisfy Theorem~\ref{thm.red 2 integrality} (but a complex construction exists), and ``?'' if the parameters satisfy Theorem~\ref{thm.red 2 integrality}, but no real construction is known.
In the end, our table identifies the parameters of every known conference graph tabulated in~\cite{Brouwer:15} (here, $k$ is simply $M-1$), meaning entries with question marks are, in fact, open problems.

We note that there are constructions of ETFs with parameters $(2,4)$, $(3,6)$ and $(5,10)$ that are not based on conference matrices.
For example, the $(2,4)$ ETFs are maximal, and there is a construction in terms of translations and modulations of a fiducial vector.
Also, a $(3,6)$ ETF can be constructed from the icosahedron: partition the 12 vertices into 6 antipodal pairs, and then take a representative from each pair.
Finally, a $(5,10)$ ETF can be obtained as a small-dimensional instance of the Tremain ETFs~\cite{FickusJMP:16}, and of the ETFs from hyperovals~\cite{FickusMJ:16}.
We did not to include these isolated constructions in the Notes column of Table~\ref{table.etfs red=2}.

\section*{Acknowledgments}

This work was supported by NSF DMS 1321779, AFOSR F4FGA05076J002 and an AFOSR Young Investigator Research Program award.
The views expressed in this article are those of the authors and do not reflect the official policy or position
of the United States Air Force, Department of Defense, or the U.S.\ Government.

\newpage
\footnotesize{
\begin{longtable}{rrrrl}
\caption{\label{table.real etfs red>2}Existence of real ETFs with $N>2M$}\\
$\quad M$& $\quad N$& $\quad \alpha$& $\quad k$& Notes\\\hline\endhead
6	&	16	&	3	&	6	&	MF(2,1), Steiner BIBD(4,2,1), Steiner AG(2,2), RSHCD(16), QSD(6,2,1)	\\
7	&	28	&	3	&	10	&	Steiner BIBD(7,3,1), Steiner PG(2,2), GQ(3,9)	\\
15	&	36	&	5	&	16	&	RSHCD(36), Tremain(7), QSD(15,3,1)	\\
19	&	76	&	5	&	32	&	DNE	\\
20	&	96	&	5	&	40	&	DNE	\\
21	&	126	&	5	&	52	&	GQ(5,25)	\\
22	&	176	&	5	&	72	&	SRG	\\
23	&	276	&	5	&	112	&	SRG	\\
28	&	64	&	7	&	30	&	MF(2,2), Steiner BIBD(8,2,1), Steiner AG(2,3), RSHCD(64), QSD(28,4,1)	\\
35	&	120	&	7	&	54	&	Steiner BIBD(15,3,1), Steiner PG(2,3)	\\
37	&	148	&	7	&	66	&		\\
41	&	246	&	7	&	108	&		\\
42	&	288	&	7	&	126	&		\\
43	&	344	&	7	&	150	&	GQ(7,49)	\\
45	&	100	&	9	&	48	&	RSHCD(100), QSD(45,5,1)	\\
45	&	540	&	7	&	234	&		\\
46	&	736	&	7	&	318	&		\\
47	&	1128	&	7	&	486	&	DNE	\\
51	&	136	&	9	&	64	&	Tremain(15)	\\
57	&	190	&	9	&	88	&		\\
61	&	244	&	9	&	112	&		\\
63	&	280	&	9	&	128	&	SRG	\\
66	&	144	&	11	&	70	&	Steiner BIBD(12,2,1), RSHCD(144), QSD(66,6,1), QSD(66,30,29)	\\
66	&	352	&	9	&	160	&		\\
69	&	460	&	9	&	208	&		\\
71	&	568	&	9	&	256	&		\\
72	&	640	&	9	&	288	&		\\
73	&	730	&	9	&	328	&	GQ(9,81)	\\
75	&	1000	&	9	&	448	&		\\
76	&	1216	&	9	&	544	&		\\
77	&	210	&	11	&	100	&	SRG	\\
88	&	320	&	11	&	150	&		\\
91	&	196	&	13	&	96	&	RSHCD(196), QSD(91,7,1)	\\
91	&	364	&	11	&	170	&		\\
99	&	540	&	11	&	250	&	Steiner BIBD(45,5,1)	\\
101	&	606	&	11	&	280	&		\\
106	&	848	&	11	&	390	&		\\
109	&	1090	&	11	&	500	&		\\
110	&	1200	&	11	&	550	&		\\
117	&	378	&	13	&	180	&		\\
120	&	256	&	15	&	126	&	MF(2,3), Steiner BIBD(16,2,1), Steiner AG(2,4), RSHCD(256), QSD(120,8,1)	\\
127	&	508	&	13	&	240	&		\\
130	&	560	&	13	&	264	&		\\
141	&	376	&	15	&	182	&		\\
141	&	846	&	13	&	396	&		\\
143	&	924	&	13	&	432	&		\\
145	&	406	&	15	&	196	&		\\
148	&	1184	&	13	&	552	&		\\
153	&	324	&	17	&	160	&	RSHCD(324)	\\
155	&	496	&	15	&	238	&	Steiner BIBD(31,3,1), Steiner PG(2,4)	\\
165	&	616	&	15	&	294	&		\\
169	&	676	&	15	&	322	&		\\
177	&	826	&	15	&	392	&		\\
183	&	976	&	15	&	462	&	Steiner BIBD(61,5,1)	\\
185	&	1036	&	15	&	490	&		\\
187	&	528	&	17	&	256	&	Tremain(31)	\\
190	&	400	&	19	&	198	&	Steiner BIBD(20,2,1), RSHCD(400)	\\
190	&	1216	&	15	&	574	&	Steiner BIBD(76,6,1)	\\
217	&	868	&	17	&	416	&		\\
221	&	936	&	17	&	448	&		\\
231	&	484	&	21	&	240	&		\\
247	&	780	&	19	&	378	&	Steiner BIBD(39,3,1)	\\
266	&	1008	&	19	&	486	&		\\
271	&	1084	&	19	&	522	&		\\
276	&	576	&	23	&	286	&	Steiner BIBD(24,2,1), RSHCD(576)	\\
276	&	736	&	21	&	360	&		\\
287	&	820	&	21	&	400	&	Tremain(39)	\\
301	&	946	&	21	&	460	&		\\
309	&	1030	&	21	&	500	&		\\
325	&	676	&	25	&	336	&	RSHCD(676)	\\
345	&	990	&	23	&	484	&		\\
365	&	876	&	25	&	432	&		\\
378	&	784	&	27	&	390	&	Steiner BIBD(28,2,1), RSHCD(784)	\\
435	&	900	&	29	&	448	&	RSHCD(900)	\\
456	&	1216	&	27	&	598	&		\\
493	&	1190	&	29	&	588	&		\\
496	&	1024	&	31	&	510	&	MF(2,4), Steiner BIBD(32,2,1), Steiner AG(2,5), RSHCD(1024), QSD(496,16,1)	\\
561	&	1156	&	33	&	576	&		\\
630	&	1296	&	35	&	646	&	Steiner BIBD(36,2,1), RSHCD(1296)	
\end{longtable}

\newpage
\begin{longtable}{rrl}
\caption{\label{table.etfs red>2}Existence of ETFs with $N>2M$}\\
$\quad M$& $\quad N$& Notes\\\hline\endhead
3	&	7	&	PG(2,2), Paley(2), SkewPaley(8), SkewHadamard(8)	\\
3	&	9	&	Steiner BIBD(3,2,1), GQ(2,4), SCHCD(9), maximal	\\
4	&	13	&	PG(2,3), Cyclo4prime(1)	\\
4	&	16	&	maximal	\\
5	&	11	&	Paley(3), SkewPaley(12)	\\
5	&	21	&	PG(2,4)	\\
5	&	25	&	maximal	\\
6	&	16	&	MF(2,1), Steiner BIBD(4,2,1), Steiner AG(2,2), RSHCD(16), SCHCD(16), QSD(6,2,1)	\\
6	&	31	&	PG(2,5)	\\
6	&	36	&	maximal	\\
7	&	15	&	PG(3,2), TPP(3), SkewHadamard(16)	\\
7	&	28	&	Steiner BIBD(7,3,1), Steiner PG(2,2), GQ(3,9) 	\\
7	&	49	&	maximal	\\
8	&	57	&	PG(2,7)	\\
8	&	64	&	maximal	\\
9	&	19	&	Paley(5), SkewPaley(20)	\\
9	&	37	&	Cyclo4(3)	\\
9	&	73	&	PG(2,8), Cyclo8(3,1)	\\
9	&	81	&	maximal	\\
10	&	25	&	Steiner BIBD(5,2,1), SCHCD(25)	\\
10	&	91	&	PG(2,9)	\\
10	&	100	&	maximal	\\
11	&	23	&	Paley(6), SkewPaley(24)	\\
11	&	121	&	maximal	\\
12	&	45	&	MF(3,1), Steiner BIBD(9,3,1), Steiner AG(3,2), Steiner Unitals(2)	\\
12	&	133	&	PG(2,11)	\\
12	&	144	&	maximal	\\
13	&	27	&	Paley(7), SkewPaley(28)	\\
13	&	40	&	PG(3,3)	\\
13	&	65	&	Steiner BIBD(13,4,1), Steiner PG(3,2), GQ(4,16)	\\
13	&	169	&	maximal	\\
14	&	183	&	PG(2,13)	\\
14	&	196	&	maximal	\\
15	&	31	&	PG(4,2), Paley(8), H(1), SkewPaley(32), SkewHadamard(32)	\\
15	&	36	&	Sp(1), Steiner BIBD(6,2,1), RSHCD(36), SCHCD(36), Tremain(7), QSD(15,3,1)	\\
15	&	225	&	maximal	\\
16	&	256	&	maximal	\\
17	&	35	&	TPP(5)	\\
17	&	273	&	PG(2,16)	\\
17	&	289	&	maximal	\\
18	&	307	&	PG(2,17)	\\
19	&	76	&	Hyperoval1(4)	\\
19	&	361	&	maximal	\\
20	&	96	&	MF(4,1), Steiner BIBD(16,4,1), Steiner AG(4,2), Hyperoval2(4)	\\
20	&	381	&	PG(2,19)	\\
21	&	43	&	Paley(11), H(2), SkewPaley(44)	\\
21	&	49	&	Steiner BIBD(7,2,1), SCHCD(49)	\\
21	&	85	&	PG(3,4)	\\
21	&	126	&	Steiner BIBD(21,5,1), Steiner PG(4,2), GQ(5,25)	\\
22	&	55	&	Tremain(9)	\\
22	&	176	&	SRG	\\
23	&	47	&	Paley(12), SkewPaley(48)	\\
23	&	276	&	SRG	\\
24	&	553	&	PG(2,23)	\\
24	&	576	&	maximal	\\
25	&	101	&	Cyclo4(5)	\\
26	&	91	&	Steiner BIBD(13,3,1)	\\
26	&	651	&	PG(2,25)	\\
28	&	64	&	MF(2,2), Steiner BIBD(8,2,1), Steiner AG(2,3), RSHCD(64), SCHCD(64), QSD(28,4,1)	\\
28	&	109	&	Cyclo4prime(5)	\\
28	&	757	&	PG(2,27)	\\
28	&	784	&	maximal	\\
29	&	59	&	Paley(15), SkewPaley(60)	\\
30	&	175	&	MF(5,1), Steiner BIBD(25,5,1), Steiner AG(5,2)	\\
30	&	871	&	PG(2,29)	\\
31	&	63	&	PG(5,2), TPP(7), SkewHadamard(64)	\\
31	&	156	&	PG(3,5)	\\
31	&	217	&	Steiner BIBD(31,6,1), Steiner PG(5,2)	\\
32	&	993	&	PG(2,31)	\\
33	&	67	&	Paley(17), SkewPaley(68)	\\
33	&	133	&	difference set	\\
33	&	1057	&	PG(2,32)	\\
35	&	71	&	Paley(18), SkewPaley(72)	\\
35	&	120	&	Steiner BIBD(15,3,1), Steiner PG(2,3)	\\
35	&	1225	&	maximal	\\
36	&	81	&	Steiner BIBD(9,2,1), SCHCD(81)	\\
39	&	79	&	Paley(20), SkewPaley(80)	\\
40	&	105	&	Tremain(13)	\\
40	&	121	&	PG(4,3)	\\
41	&	83	&	Paley(21), SkewPaley(84)	\\
43	&	344	&	GQ(7,49)	\\
45	&	100	&	Steiner BIBD(10,2,1), RSHCD(100), SCHCD(100), QSD(45,5,1)	\\
49	&	99	&	TPP(9)	\\
49	&	197	&	Cyclo4(7)	\\
50	&	225	&	Steiner BIBD(25,4,1)	\\
51	&	103	&	Paley(26), SkewPaley(104)	\\
51	&	136	&	Tremain(15)	\\
53	&	107	&	Paley(27), SkewPaley(108)	\\
55	&	121	&	Steiner BIBD(11,2,1), SCHCD(121)	\\
56	&	441	&	MF(7,1), Steiner AG(7,2)	\\
57	&	190	&	Steiner BIBD(19,3,1)	\\
57	&	400	&	PG(3,7)	\\
57	&	513	&	Steiner PG(7,2), GQ(8,64)	\\
63	&	127	&	PG(6,2), Paley(32), H(5), SkewPaley(128), SkewHadamard(128)	\\
63	&	280	&	Steiner BIBD(28,4,1), Steiner Unitals(3), Steiner Denniston(2,3), SRG	\\
65	&	131	&	Paley(33), SkewPaley(132)	\\
66	&	144	&	difference set, Steiner BIBD(12,2,1), RSHCD(144), SCHCD(144), QSD(66,6,1), QSD(66,30,29)	\\
69	&	139	&	Paley(35), SkewPaley(140)	\\
70	&	231	&	Steiner BIBD(21,3,1)	\\
71	&	143	&	TPP(11)	\\
71	&	568	&	Hyperoval1(8)	\\
72	&	640	&	MF(8,1), Steiner AG(8,2), Hyperoval2(8)	\\
73	&	585	&	PG(3,8)	\\
73	&	730	&	Steiner PG(8,2), GQ(9,81)	\\
75	&	151	&	Paley(38), SkewPaley(152)	\\
77	&	210	&	Tremain(19), SRG	\\
78	&	169	&	Steiner BIBD(13,2,1), SCHCD(169)	\\
81	&	163	&	Paley(41), SkewPaley(164)	\\
82	&	451	&	Steiner BIBD(41,5,1)	\\
83	&	167	&	Paley(42), SkewPaley(168)	\\
85	&	341	&	PG(4,4)	\\
88	&	320	&	difference set	\\
89	&	179	&	Paley(45), SkewPaley(180)	\\
90	&	891	&	MF(9,1), Steiner AG(9,2)	\\
91	&	196	&	Steiner BIBD(14,2,1), RSHCD(196), SCHCD(196), QSD(91,7,1)	\\
91	&	820	&	PG(3,9)	\\
91	&	1001	&	Steiner PG(9,2)	\\
92	&	253	&	Tremain(21)	\\
95	&	191	&	Paley(48), SkewPaley(192)	\\
99	&	199	&	Paley(50), SkewPaley(200)	\\
99	&	540	&	Steiner BIBD(45,5,1)	\\
100	&	325	&	Steiner BIBD(25,3,1)	\\
105	&	211	&	Paley(53), SkewPaley(212)	\\
105	&	225	&	Steiner BIBD(15,2,1), SCHCD(225)	\\
111	&	223	&	Paley(56), H(7), SkewPaley(224)	\\
111	&	481	&	Steiner BIBD(37,4,1)	\\
113	&	227	&	Paley(57), SkewPaley(228)	\\
117	&	378	&	MF(3,2), Steiner BIBD(27,3,1), Steiner AG(3,3)	\\
119	&	239	&	Paley(60), SkewPaley(240)	\\
120	&	256	&	MF(2,3), Steiner BIBD(16,2,1), Steiner AG(2,4), RSHCD(256), SCHCD(256), QSD(120,8,1)	\\
121	&	243	&	Paley(61), SkewPaley(244)	\\
121	&	364	&	PG(5,3)	\\
125	&	251	&	Paley(63), SkewPaley(252)	\\
126	&	351	&	Sp(2), Tremain(25)	\\
127	&	255	&	PG(7,2), SkewHadamard(256)	\\
130	&	560	&	Steiner BIBD(40,4,1), Steiner PG(3,3)	\\
131	&	263	&	Paley(66), SkewPaley(264)	\\
135	&	271	&	Paley(68), SkewPaley(272)	\\
136	&	289	&	Steiner BIBD(17,2,1), SCHCD(289)	\\
141	&	283	&	Paley(71), H(8), SkewPaley(284)	\\
143	&	924	&	Steiner BIBD(66,6,1)	\\
145	&	406	&	Tremain(27)	\\
153	&	307	&	Paley(77), SkewPaley(308)	\\
153	&	324	&	difference set, Steiner BIBD(18,2,1), RSHCD(324), SCHCD(324)	\\
155	&	311	&	Paley(78), SkewPaley(312)	\\
155	&	496	&	Steiner BIBD(31,3,1), Steiner PG(2,4)	\\
156	&	781	&	PG(4,5)	\\
161	&	323	&	TPP(17)	\\
165	&	331	&	Paley(83), SkewPaley(332)	\\
169	&	677	&	Cyclo4(13)	\\
171	&	343	&	Paley(86), SkewPaley(344)	\\
171	&	361	&	Steiner BIBD(19,2,1), SCHCD(361)	\\
173	&	347	&	Paley(87), SkewPaley(348)	\\
176	&	561	&	Steiner BIBD(33,3,1)	\\
179	&	359	&	Paley(90), SkewPaley(360)	\\
183	&	367	&	Paley(92), SkewPaley(368)	\\
183	&	976	&	Steiner BIBD(61,5,1)	\\
187	&	528	&	Tremain(31)	\\
189	&	379	&	Paley(95), SkewPaley(380)	\\
190	&	400	&	Steiner BIBD(20,2,1), RSHCD(400), SCHCD(400)	\\
190	&	1216	&	Steiner BIBD(76,6,1)	\\
191	&	383	&	Paley(96), SkewPaley(384)	\\
196	&	833	&	Steiner BIBD(49,4,1)	\\
208	&	1105	&	Steiner BIBD(65,5,1), Steiner Unitals(4)	\\
209	&	419	&	Paley(105), SkewPaley(420)	\\
210	&	441	&	Steiner BIBD(21,2,1), SCHCD(441)	\\
210	&	595	&	Tremain(33)	\\
215	&	431	&	Paley(108), SkewPaley(432)	\\
219	&	439	&	Paley(110), SkewPaley(440)	\\
221	&	443	&	Paley(111), SkewPaley(444)	\\
221	&	936	&	Steiner BIBD(52,4,1), Steiner Denniston(2,4)	\\
222	&	703	&	Steiner BIBD(37,3,1)	\\
225	&	901	&	difference set	\\
231	&	463	&	Paley(116), SkewPaley(464)	\\
231	&	484	&	Steiner BIBD(22,2,1), SCHCD(484)	\\
233	&	467	&	Paley(117), SkewPaley(468)	\\
239	&	479	&	Paley(120), SkewPaley(480)	\\
243	&	487	&	Paley(122), SkewPaley(488)	\\
245	&	491	&	Paley(123), SkewPaley(492)	\\
247	&	780	&	Steiner BIBD(39,3,1)	\\
249	&	499	&	Paley(125), SkewPaley(500)	\\
251	&	503	&	Paley(126), SkewPaley(504)	\\
253	&	529	&	Steiner BIBD(23,2,1), SCHCD(529)	\\
255	&	511	&	PG(8,2), SkewHadamard(512)	\\
260	&	741	&	Tremain(37)	\\
261	&	523	&	Paley(131), SkewPaley(524)	\\
273	&	547	&	Paley(137), SkewPaley(548)	\\
276	&	576	&	difference set, Steiner BIBD(24,2,1), RSHCD(576), SCHCD(576)	\\
281	&	563	&	Paley(141), SkewPaley(564)	\\
285	&	571	&	Paley(143), SkewPaley(572)	\\
287	&	575	&	TPP(23)	\\
287	&	820	&	Tremain(39)	\\
293	&	587	&	Paley(147), SkewPaley(588)	\\
299	&	599	&	Paley(150), SkewPaley(600)	\\
300	&	625	&	Steiner BIBD(25,2,1), SCHCD(625)
\end{longtable}

\newpage
\begin{longtable}{rrcl}
\caption{\label{table.etfs red=2}Existence of ETFs with $N=2M$}\\
$\quad M$& $\quad N$& \quad $\mathbb{R}?$ \quad & Notes\\\hline \endhead
2	&	4	&	-	&	(f), (g)	\\
3	&	6	&	+	&	(a)	\\
4	&	8	&	-	&	(f), (g)	\\
5	&	10	&	+	&	(a), (d)	\\
6	&	12	&	-	&	(f)	\\
7	&	14	&	+	&	(a)	\\
8	&	16	&	-	&	(g)	\\
9	&	18	&	+	&	(a)	\\
10	&	20	&	-	&	(f)	\\
12	&	24	&	-	&	(f)	\\
13	&	26	&	+	&	(a), (e)	\\
14	&	28	&	-	&	(f)	\\
15	&	30	&	+	&	(a)	\\
16	&	32	&	-	&	(f), (g)	\\
19	&	38	&	+	&	(a)	\\
21	&	42	&	+	&	(a)	\\
22	&	44	&	-	&	(f)	\\
23	&	46	&	+	&	(b), (c)	\\
24	&	48	&	-	&	(f)	\\
25	&	50	&	+	&	(a), (d)	\\
27	&	54	&	+	&	(a)	\\
30	&	60	&	-	&	(f)	\\
31	&	62	&	+	&	(a)	\\
32	&	64	&	-	&	(g)	\\
33	&	66	&	?	&		\\
34	&	68	&	-	&	(f)	\\
36	&	72	&	-	&	(f)	\\
37	&	74	&	+	&	(a)	\\
40	&	80	&	-	&	(f)	\\
41	&	82	&	+	&	(a), (d), (e)	\\
42	&	84	&	-	&	(f)	\\
43	&	86	&	?	&		\\
45	&	90	&	+	&	(a)	\\
49	&	98	&	+	&	(a)	\\
51	&	102	&	+	&	(a)	\\
52	&	104	&	-	&	(f)	\\
54	&	108	&	-	&	(f)	\\
55	&	110	&	+	&	(a)	\\
57	&	114	&	+	&	(a)	\\
59	&	118	&	?	&		\\
61	&	122	&	+	&	(a)	\\
63	&	126	&	+	&	(a), (e)	\\
64	&	128	&	-	&	(f), (g)	\\
66	&	132	&	-	&	(f)	\\
69	&	138	&	+	&	(a)	\\
70	&	140	&	-	&	(f)	\\
73	&	146	&	?	&		\\
75	&	150	&	+	&	(a)	\\
76	&	152	&	-	&	(f)	\\
77	&	154	&	?	&		\\
79	&	158	&	+	&	(a)	\\
82	&	164	&	-	&	(f)	\\
84	&	168	&	-	&	(f)	\\
85	&	170	&	+	&	(a), (e)	\\
87	&	174	&	+	&	(a)	\\
90	&	180	&	-	&	(f)	\\
91	&	182	&	+	&	(a)	\\
93	&	186	&	?	&		\\
96	&	192	&	-	&	(f)	\\
97	&	194	&	+	&	(a)	\\
99	&	198	&	+	&	(a)	\\
100	&	200	&	-	&	(f)	\\
103	&	206	&	?	&		\\
106	&	212	&	-	&	(f)	\\
111	&	222	&	?	&		\\
112	&	224	&	-	&	(f)	\\
113	&	226	&	+	&	(d)	\\
114	&	228	&	-	&	(f)	\\
115	&	230	&	+	&	(a)	\\
117	&	234	&	+	&	(a)	\\
120	&	240	&	-	&	(f)	\\
121	&	242	&	+	&	(a)	\\
122	&	244	&	-	&	(f)	\\
123	&	246	&	?	&		\\
126	&	252	&	-	&	(f)	\\
128	&	256	&	-	&	(g)	\\
129	&	258	&	+	&	(a)	\\
131	&	262	&	?	&		\\
132	&	264	&	-	&	(f)	\\
133	&	266	&	?	&		\\
135	&	270	&	+	&	(a)	\\
136	&	272	&	-	&	(f)	\\
139	&	278	&	+	&	(a)	\\
141	&	282	&	+	&	(a)	\\
142	&	284	&	-	&	(f)	\\
145	&	290	&	+	&	(a), (e)	\\
147	&	294	&	+	&	(a)
\end{longtable}

}

\end{document}